\documentclass[12pt,a4paper]{article}
\usepackage{amssymb}
\usepackage{amsfonts}
\usepackage{amsmath}
\usepackage{amsthm}
\usepackage{enumerate}
\usepackage{srcltx}
\newtheorem{theorem}{Theorem}[section]
\newtheorem{proposition}{Proposition}[section]
\newtheorem{lemma}{Lemma}[section]

\newtheorem{corollary}{Corollary}[section]

\newtheorem{remark}{Remark}[section]
\newcommand{\wh}{\widehat}
\newcommand{\wt}{\widetilde}

\newcommand\cF{{\cal F}}
\newcommand\cL{{\cal L}}

\newcommand\cK{{\cal K}}

\newcommand\cV{{\cal V}}

\def\bbr{{\mathbb R}}
\def\x{{\bf x}}

\newcommand{\aO}{\mbox{O}}

\def\text#1{\hbox{#1}}
\def\proof{{\noindent \bf Proof. }}
\def\endproof{\mbox{\ $\qed$}}

\def\E{{\bf E}}
\def\P{{\bf P}}
\def\C{{\bf C}}

\def\I{{\bf I}}

\def\k{{\bf k}}

\def\A{{\bf A}}
\def\c{{\bf c}}
\def\B{{\bf B}}

\def\Chi{{\bf 1}}
\def\d{\mathrm{d}}
\def\build #1_#2{\mathrel{\mathop{\kern 0pt #1}\limits_\zs{#2}}}
\newcommand{\zs}[1]{{\mathchoice{#1}{#1}{\lower.25ex\hbox{$\scriptstyle#1$}}
{\lower0.25ex\hbox{$\scriptscriptstyle#1$}}}}

\numberwithin{equation}{section}
\def\proof{{\noindent \bf Proof. }}
\def\endproof{\mbox{\ $\qed$}}
\begin{document}
\title{
Uniform concentration inequality for ergodic
 diffusion processes observed at discrete times. 
\thanks{The second author is partially supported by the RFFI-Grant  09-01-00172-a.}
}

\author{L. Galtchouk
\thanks{
Department of Mathematics,
 Strasbourg University
7, rue Rene Descartes,
67084, Strasbourg, France, 
e-mail: leonid.galtchouk@math.unistra.fr }
 \and 
S. Pergamenshchikov
\thanks{
 Laboratoire de Math\'ematiques Raphael Salem,
 Avenue de l'Universit\'e, BP. 12,            
 Universit\'e de Rouen,                  
F76801, Saint Etienne du Rouvray, Cedex France, 
e-mail: Serge.Pergamenchtchikov@univ-rouen.fr}
}

\maketitle

\begin{abstract}
In this paper a
concentration inequality is
 proved for the deviation in the ergodic theorem in the case 
of discrete time observations of diffusion processes. 
The proof is based on the geometric ergodicity
property for diffusion processes. As an application 
we consider the nonparametric pointwise estimation problem 
for the drift coefficient under discrete time observations.
\end{abstract}


\par


\bibliographystyle{plain}

\newpage

\section{Introduction}\label{sec:In}

We consider the process $(y_\zs{t})_\zs{t\ge 0}$ governed by the
stochastic differential equation
\begin{equation}\label{sec:In.1}
\d y_\zs{t}=S(y_\zs{t})\,\d t +\sigma(y_\zs{t})\d W_\zs{t}\,,
\quad 0\le t\le T\,,
\end{equation}
where $(W_\zs{t}, \cF_\zs{t})_\zs{t\ge 0}$ is a standard Wiener process, 
$y_\zs{0}$ is a
initial condition and $\vartheta=(S,\sigma)$ are unknown functions.
For this model we consider the pointwise estimation problem  for the function 
$S$ at a fixed point 
$x_\zs{0}\in\bbr$ (i.e. $S(x_\zs{0})$), on the basis of the discrete time observations of 
the process \eqref{sec:In.1}, i.e.  
\begin{equation}\label{sec:In.2}
(y_\zs{t_\zs{j}})_\zs{1\le j\le N}\,,
\end{equation}
where $t_\zs{j}=j\delta$, $N=[T/\delta]$ and $\delta$ is some positive fixed 
observation
frequency  which will be specified later.
Usually, for this problem one uses 
 kernel estimators $\wh{S}_\zs{N}(x_\zs{0})$ defined as
\begin{equation}\label{sec:In.3}
\wh{S}_\zs{N}(x_\zs{0})=\frac{\sum_\zs{k=1}^{N}\,
\psi_\zs{h,x_\zs{0}}(y_\zs{t_\zs{k}})
\,\Delta y_\zs{t_\zs{k}}}
{\sum_\zs{k=1}^{N}\,\psi_\zs{h,x_\zs{0}}(y_\zs{t_\zs{k}})\,  \Delta t_\zs{k}}\,,
\quad 
\psi_\zs{h,x_\zs{0}}(y)=\frac{1}{h}\,\Psi\left(\frac{y-x_\zs{0}}{h}\right)\,,
\end{equation}
where $\Psi(y)$ is a kernel function which equals to zero for $|y|\ge 2$ and
 will be specified later,
 $0<h<1$ is a bandwidth, 
$\Delta y_\zs{t_\zs{k}}=y_\zs{t_\zs{k}}-y_\zs{t_\zs{k-1}}$
and $\Delta t_\zs{k}=\delta$.

Main difficulty in this estimator is that the denominator is random. Therefore, 
to obtain the convergence rate for this estimator we have to study the behavior
of the denominator, more precisely, one needs to show that
$$
\sum_\zs{k=1}^{N}\,\psi_\zs{h,x_\zs{0}}(y_\zs{t_\zs{k}})\Delta t_\zs{k}
\approx \pi_\zs{\vartheta}(\psi_\zs{h,x_\zs{0}}) h T
\quad\mbox{as}\quad T\to\infty
\,,
$$
where
\begin{equation}\label{sec:In.4}
\pi_\zs{\vartheta}(\psi_\zs{h,x_\zs{0}})=
\int_\zs{\bbr}\,\psi_\zs{h,x_\zs{0}}(y)\,q_\zs{\vartheta}(y)\,\d y
\end{equation}
and $q_\zs{\vartheta}$ is the ergodic density defined in 
\eqref{sec:Mr.2}.

Unfortunately, the ergodic theorem does not permit to obtain this kind of result because the times
$t_\zs{k}$ and the bandwidth $h$ depend on $T$.
Usually one obtains such properties through concentration inequalities for the deviation
in the ergodic theorem, i.e. one needs to study the limit behavior
of the deviation
\begin{equation}\label{sec:In.5}
D_\zs{T}(\phi)
=
\sum_\zs{k=1}^{N}\,
\left(
\phi(y_\zs{t_\zs{k}})
-
\pi_\zs{\vartheta}(\phi)
\right)
\,\Delta t_\zs{k}
\end{equation}
for some functions $\phi$ which can be dependent on $T$, for example, 
$\phi(\cdot)=\psi_\zs{h,x_\zs{0}}(\cdot)$.
More precisely, we need to
show, that for any $\varepsilon>0$ and for any $m>0$,
uniformly over $\vartheta$,
\begin{equation}\label{sec:In.6}
\lim_\zs{T\to\infty}\,
T^{m}\,
\P_\zs{\vartheta}
\left(
|
D_\zs{T}(\psi_\zs{h,x_\zs{0}})
|
> \varepsilon T
\right)=0\,,
\end{equation}
where $\P_\zs{\vartheta}$ is the law of the process 
$(y_\zs{t})_\zs{t\ge 0}$ under the coefficients $\vartheta=(S,\sigma)$.
Usually, to get  properties of type \eqref{sec:In.6} 
one needs to establish 
an exponential inequality for the deviations  
\eqref{sec:In.5}.

There are a number of papers devoted to concentration inequalities
for functions of independent random variables (we refer the reader
to \cite{BoMa} and references therein),
for functions of dependent random variables (see \cite{DeDo},
\cite{DePr}, \cite{Ma}). For Markov chains such inequalities were obtained
in \cite{BeCl}. For continuous time Markov processes an exponential
concentration inequality was obtained in \cite{CaGu} (see also references therein). Some applications of concentration inequalities to statistics are presented in \cite{Mas}. 
Concentration inequalities for diffusion processes are given in
\cite{GaPe2}, \cite{Pe}, \cite{Ve}.

 For statistical applications, we need uniform upper bounds for the tail distribution 
over functions $\phi$ like to the exponential bounds in \cite{GaPe2}.
We can not apply directly the method  from \cite{GaPe2}, since there
it is based on the continuous times version of the Ito 
formula. In this paper we apply this approach through 
uniform (over the functions $S$)
geometric ergodicity. We recall (see \cite{MeTw}), that the geometric ergodicity
yields a geometric rate in the convergence
$$
\lim_\zs{t\to\infty}\,\E_\zs{\vartheta}(g(y_\zs{t})|y_\zs{0}=x)=
\pi_\zs{\vartheta}(g)
$$
for any integrable functions $g$ and any initial value $x\in\bbr$.
Here $\E_\zs{\vartheta}$ denotes the expectation with respect to the distribution $\P_\zs{\vartheta}$. 
In \cite{GaPe4} through the Lyapunov functions method 
it is shown that the process \eqref{sec:In.1} is geometrically ergodic
uniformly over functions $\vartheta=(S,\sigma)$ from the functional class $\Theta$ 
defined  in
\eqref{sec:Mr.1}.

The paper is organized as follows. In the next section 
we formulate the main results. In Section~\ref{sec:Pa} 
we introduce all the necessary parameters. In Section~
\ref{sec:Ci} we show a concentration inequality in ergodic theorem
for the continuous observations of the process \eqref{sec:In.1}.
In Section~\ref{sec:Ge} we announce the uniform geometric ergodic property
for the process  \eqref{sec:In.1}.
In Section~\ref{Bi}
we give the Burkh\"older inequality
for dependent random variables.
In Section~\ref{sec:Pr} we prove all main results.
The Appendix contains the proofs of some auxiliary
results.

\section{Main results}\label{sec:Mr}

First we describe the functional class $\Theta$ for functions $\vartheta=(S,\sigma)$
defined in \cite{GaPe4}. We start with some real numbers 
$\x_\zs{*}\ge 1$, $M>0$ and  $L>1$ for
which we denote by $\Sigma_\zs{L,M}$
the class of functions $S$ from $\C^{1}(\bbr)$ such that
$$
\sup_\zs{|x|\le \x_\zs{*}}
\left(|S(x)|+|\dot{S}(x)|\right)\le M
$$
and
$$
-L\le \inf_\zs{|x|\ge \x_\zs{*}}\dot{S}(x)\le \sup_\zs{|x|\ge \x_\zs{*}}\dot{S}(x)
\le -L^{-1}
\,.
$$
Furthermore, for some fixed numbers
 $0<\sigma_\zs{\min}\le \sigma_\zs{\max}<\infty$, we denote by $\cV$
the class of the functions $\sigma$ from $\C^{2}(\bbr)$ such that
\begin{align*}
\sigma_\zs{min}&\le\,\inf_\zs{x\in\bbr}
\min\left(|\sigma(x)|\,,\,|\dot{\sigma}(x)|\,,\,|\ddot{\sigma}(x)|\right)
\\[2mm] 
&
\le \sup_\zs{x\in\bbr}
\max\left(|\sigma(x)|\,,\,|\dot{\sigma}(x)|\,,\,|\ddot{\sigma}(x)|\right)
\le\sigma_\zs{max}\,.
\end{align*}

\noindent Finally, we set
\begin{equation}\label{sec:Mr.1}
\Theta=\Sigma_\zs{L,M}\times \cV\,.
\end{equation}

\noindent 
It should be noted (see, for example, \cite{GiSk}), that
 for any $\vartheta=(S,\sigma)\in\Theta$, 
 the equation \eqref{sec:In.1} has a unique strong  solution 
which is a ergodic process with the invariant density 
$q_\zs{\vartheta}$ defined as
\begin{equation}\label{sec:Mr.2}
q_\zs{\vartheta}(x)=
\left(\int_\zs{\bbr}\sigma^{-2}(z)\,
e^{\wt{S}(z)}\d z\right)^{-1}\,
\sigma^{-2}(x)\,e^{\wt{S}(x)}\,,
 \end{equation}
where $\wt{S}(x)=2\int_\zs{0}^{x}\,S_\zs{1}(v)\d v$
and  $S_\zs{1}(x)=S(x)/\sigma^{2}(x)$.

Now we describe the functional classes for the functions $\phi$. First,
for any parameters $\nu_\zs{0}>0$ and $\nu_\zs{1}>0$ we set
\begin{equation}\label{sec:Mr.2-1}
\cV_\zs{\nu_\zs{0},\nu_\zs{1}}=
\left\{\phi\in\C(\bbr)\,:\,|\phi|_\zs{1}\le \nu_\zs{0}\,,\,
|\phi|_\zs{*}\le \nu_\zs{1}
\right\}\,,
\end{equation}
where $|\phi|_\zs{1}=\int_\zs{\bbr}\,|\phi(y)|\,\d y$ and 
$|\phi|_\zs{*}=\sup_\zs{y\in\bbr}|\phi(y)|$.

For any 
function $\phi$ from $\C^{2}(\bbr)$ 
we denote by $\cL_\zs{\vartheta}(\phi)$ the generator operator for the process
\eqref{sec:In.1}, i.e.
$$
\cL_\zs{\vartheta}(\phi)(y)=S(y)\dot{\phi}(y)+\frac{\sigma^{2}(y)}{2}\ddot{\phi}(y)
\,.
$$
\noindent 
Using this notation, we set
\begin{equation}\label{sec:Pa.1}
\mu(\phi)=\sup_\zs{\vartheta\in\Theta}\|\cL_\zs{\vartheta}(\phi)\|_\zs{*}
\,
\quad\mbox{and}\quad
\wt{\mu}(\phi)=\sup_\zs{\vartheta\in\Theta}
\,
|\wt{\pi}_\zs{\vartheta}(\phi)|
\,,
\end{equation}
where $\wt{\pi}_\zs{\vartheta}(\phi)=\pi_\zs{\vartheta}(\cL_\zs{\vartheta}(\phi))$. 
Now for any vector 
$\nu=(\nu_\zs{0},\nu_\zs{1},\nu_\zs{2},\nu_\zs{3},\nu_\zs{4})$
from $\bbr_\zs{+}^{5}$ we set

\begin{equation}\label{sec:Mr.2-2}
\cK_\zs{\nu}=
\left\{\phi\in \cV_\zs{\nu_\zs{0},\nu_\zs{1}}\,:\,
\|\dot{\phi}\|_\zs{*}\le \nu_\zs{2}\,,\,\mu(\phi)\le \nu_\zs{3}\,,\,
\wt{\mu}(\phi)\le \nu_\zs{4}
\right\}\,.
\end{equation}

\medskip

\begin{theorem}\label{Th.sec:Mr.1}
For any vector 
$\nu=(\nu_\zs{0},\nu_\zs{1},\nu_\zs{2},\nu_\zs{3},\nu_\zs{4})$
from $\bbr_\zs{+}^{5}$
and any\\ $0<\delta\le 1$ there exist positive parameters
$z_\zs{0}=z_\zs{0}(\delta,\nu)$, $\gamma=\gamma(\delta,\nu)$ 
and $\varkappa=\varkappa(\delta,\nu)$ such that 
\begin{equation}\label{sec:Mr.3}
\sup_\zs{T\ge 1}\,
\sup_\zs{z\ge z_\zs{0}}\,
\sup_\zs{\phi\in \cK_\zs{\nu}}
\sup_\zs{\vartheta\in\Theta}\,
e^{ z\min(\varkappa z\,,\,\gamma)}
\,
\P_\zs{\vartheta}\left(
|D_\zs{T}(\phi)|\,
\ge z \sqrt{N}
\right)
\le 4
\,,
\end{equation}
where  the parameters $z_\zs{0}$, $\gamma$ and $\varkappa$ are defined 
in  \eqref{sec:Pa.6}--\eqref{sec:Pa.7}.
\end{theorem}

\noindent Now we apply this theorem to the pointwise estimation problem, i.e.
for the functions $\psi_\zs{h,x_\zs{0}}$ defined in \eqref{sec:In.3}. 
 To this end
we assume that the frequency $\delta$ in the observations
 \eqref{sec:In.2} is of the following form
\begin{equation}\label{sec:Mr.4}
\delta=\delta_\zs{T}=\frac{1}{T\, l_\zs{T}}\,,
\end{equation}
where the function $l_\zs{T}$ is such that for any $m>0$
\begin{equation}\label{sec:Mr.5}
\lim_\zs{T\to\infty}\frac{l_\zs{T}}{T^{m}}=0
\quad\mbox{and}\quad
\lim_\zs{T\to\infty}\frac{l_\zs{T}}{\ln T}=+\infty\,.
\end{equation}

\noindent 
Further, let $\epsilon=\epsilon_\zs{T}$ be a positive function satisfying the following properties
\begin{equation}\label{sec:Mr.6}
\lim_\zs{T\to\infty}\,\epsilon_\zs{T}=0, \quad
\lim_\zs{T\to\infty}\,\frac{l_\zs{T}}{T\epsilon_\zs{T}}=0
\quad\mbox{and}\quad
\lim_\zs{T\to\infty}\frac{\epsilon^{5}_\zs{T}l_\zs{T}}{\ln T}=+\infty\,.
\end{equation}
We can take, for example, for some $\iota>0$
$$
l_\zs{T}=\ln^{1+6\iota}(T+1)
\quad\mbox{and}\quad
\epsilon_\zs{T}=\frac{1}{ \ln^{\iota}(T+1)}\,.
$$
\begin{theorem}\label{Th.sec:Mr.2}
Assume that the  kernel function $\Psi$ in \eqref{sec:In.3} 
is two continuously differentiable.
Moreover, assume that the functions $\delta_\zs{T}$
and $l_\zs{T}$ satisfy the properties \eqref{sec:Mr.4} and
\eqref{sec:Mr.6}. Then there exist 
coefficients $z^{*}_\zs{0}=z^{*}_\zs{0}(\Psi)>0 $
and
$\gamma^{*}=\gamma^{*}(\Psi)>0$
 such that 
\begin{equation}\label{sec:Mr.7}
\limsup_\zs{T\to\infty}\,e^{a\gamma^{*}\,l_\zs{T}}
\,\sup_\zs{a\ge a_\zs{*}}
\sup_\zs{h\ge T^{-1/2}}
\,
\sup_\zs{\vartheta\in\Theta}
\,
\P_\zs{\vartheta}\left(
|D_\zs{T}(\psi_\zs{h,x_\zs{0}})|\,
\ge
\,a\,T
\right)
\le 4\,,
\end{equation}
where $a_\zs{*}=z^{*}_\zs{0}/l_\zs{T}$, the parameters 
$z^{*}_\zs{0}$ and $\gamma^{*}$ are given in Section~\ref{sec:Pa}.
\end{theorem}
\noindent 
This theorem implies immediately the following
\begin{corollary}\label{Co.sec:App.1}
Assume, that all conditions 
of Theorem~\ref{Th.sec:Mr.2} hold. Then,
for any $m>0$,
$$
\limsup_\zs{T\to\infty}\,
T^{m}
\,\sup_\zs{a\ge a_\zs{*}}
\sup_\zs{h\ge T^{-1/2}}
\,
\sup_\zs{\vartheta\in\Theta}
\,
\P_\zs{\vartheta}\left(
|D_\zs{T}(\psi_\zs{h,x_\zs{0}})|\,
\ge
\,a\,T
\right)
=0\,.
$$
\end{corollary}

\noindent
Now we study the deviation \eqref{sec:In.5} for the function
\begin{equation}\label{sec:Mr.7-1}
\chi_\zs{h,x_\zs{0}}(y)=
\frac{1}{h}\,\chi\left(\frac{y-x_\zs{0}}{h}\right)\,,
\end{equation}
where $\chi(y)=\Chi_\zs{\{|y|\le 1\}}$. 

\begin{theorem}\label{Th.sec:Mr.3}
Assume that the parameter $\delta$ has the form
\eqref{sec:Mr.4}. 
 Then, for any $m>0$, and for any function
 $\epsilon_\zs{T}$, satisfying the condtions
 \eqref{sec:Mr.5} and
\eqref{sec:Mr.6}
\begin{equation}\label{sec:Mr.8}
\lim_\zs{T\to\infty}\,T^{m}\,
\sup_\zs{h\ge T^{-1/2}}
\quad
\sup_\zs{\vartheta\in\Theta}
\,
\P_\zs{\vartheta}\left(
|D_\zs{T}(\chi_\zs{h,x_\zs{0}})|\,
\ge
\epsilon_\zs{T}\,T
\right)=0\,.
\end{equation}
\end{theorem}

\begin{remark}\label{Re.sec:App.1}
It is well known that to obtain the optimal rate in the estimation problem
for a differentiable function $S$ in the process \eqref{sec:In.1} one needs to
choose the bandwidth $h$ as
$$
h=T^{-1/(2\alpha+1)}
$$
with the regularity parameter $\alpha\ge 1$. 
This means that, really for the pointwise estimation problem, $h\ge T^{-1/3}$.
But in the  quadratic risk one needs to choose the parameter
$h$ as $h=T^{-1/2}$ (see \cite{GaPe}-\cite{GaPe1},\cite{GaPe3}).
\end{remark}

\medskip

\section{Parameters}\label{sec:Pa}

In this section we introduce all necessary constants and parameters. 
First, we set
\begin{equation}\label{sec:Pa.2}
\upsilon_\zs{1}=e^{\beta^{2}_\zs{1}/(4\beta_\zs{2})}
\quad\mbox{and}\quad
\upsilon_\zs{2}=
\sqrt{\pi/\beta_\zs{2}}\,
e^{\beta^{2}_\zs{1}/(4\beta_\zs{2})}\,,
\end{equation}
where
$\beta_\zs{1}=2M/\sigma^{2}_\zs{\min}$
and
$\beta_\zs{2}=1/L\sigma^{2}_\zs{\max}$. 
 Moreover, as we will see in Appendix, the
ergodic density  \eqref{sec:Mr.2} is 
uniformly bounded by $q^{*}$, where
\begin{equation}\label{sec:Pa.3}
q^{*}=
\frac{\sigma^{2}_\zs{max}}{\sigma^{2}_\zs{min}}
\,e^{\beta_\zs{1}\x_\zs{*}+\beta^{2}_\zs{1}/(4\beta_\zs{2})}
\,.
\end{equation}
\noindent
Now we set

\begin{equation}\label{sec:Pa.4}
r=r(\nu_\zs{0})=
\frac{2 \nu_\zs{0}}{\sigma^{2}_\zs{\min}}
\,
\left(1+\upsilon_\zs{1}+q^{*}
\left(\x_\zs{*}+
\upsilon_\zs{2}\right)
\right)
\,
e^{\x_\zs{*}\beta_\zs{1}}
\,,
\end{equation}
where the parameter $\nu_\zs{0}$ is defined in \eqref{sec:Mr.2-1}.
Now using this function we set
\begin{equation}\label{sec:Pa.5}
\kappa_\zs{0}=\kappa_\zs{0}(\nu_\zs{0})=
\frac{1}{108\,r^{2}(3\rho^{2}+y^{2}_\zs{0}+2\sigma^{2}_\zs{\max})}
\end{equation}
where $\rho=
\max\left(|y_\zs{0}|\,,\,\sigma_\zs{\max}\sqrt{L}\,,\,2(\x_\zs{*}+M L)\right)$.

Now for any 
$\delta>0$ and any parameter vector 
$\nu=(\nu_\zs{0},\nu_\zs{1},\nu_\zs{2},\nu_\zs{3},\nu_\zs{4})$ from
$\bbr^{5}_\zs{+}$  we set
\begin{align}\nonumber
z_\zs{0}&=
z_\zs{0}(\delta,\nu)=\delta^{3/2}\,
\max\left(2c^{*}_\zs{1}\nu_\zs{3}\,,\,
2c^{*}_\zs{2}\nu_\zs{2}\,,\, \nu_\zs{4}T^{1/2}\,,\,
\nu_\zs{1} T^{-1/2}
\right)\,,\\[5mm]\label{sec:Pa.6}
\tau&=\tau(\delta,\nu)=\delta^{3/2}\,
\max\left(c^{*}_\zs{1}\nu_\zs{3}\,,\,
c^{*}_\zs{2}\nu_\zs{2}\right)\,,
\end{align}
where 
$$
c^{*}_\zs{1}=2e^{\kappa+1}\sqrt{\frac{R(1+\rho)}{\kappa}}
\quad\mbox{and}\quad
c^{*}_\zs{2}=\sqrt{2}e\sigma_\zs{max}
\,.
$$
The parameters $R$ and $\kappa$ are 
defined in Theorem~\ref{Th.sec:Ci.0}.
Finally we set

\begin{equation}\label{sec:Pa.7}
\gamma=\frac{1}{4\tau}
\quad\mbox{and}\quad
\varkappa=\varkappa(\delta,\nu)=
\frac{9\kappa_\zs{0}(1-\delta)}{64\delta}\,.
\end{equation}

\noindent 
Now we set
\begin{equation}\label{sec:Pa.8}
M_\zs{1}=
M+L\left(\x_\zs{*}+
|x_\zs{0}|+2
\right)
\,.
\end{equation}

\noindent Now for any inegrated two times continuously differentiable
$\bbr\to\bbr$ function $\Psi$ we define 
\begin{equation}\label{sec:Pa.9}
\k_\zs{*}(\Psi)=
  \max\left(
|\dot{\Psi}|_\zs{1}
\,,\,
|\ddot{\Psi}|_\zs{1}
\,,\,
\|\Psi\|_\zs{*}
\,,\,
\|\dot{\Psi}\|_\zs{*}
\,,\,
\|\ddot{\Psi}\|_\zs{*}
\right)\,.
\end{equation}
\noindent
Using this operator we define the parameters
\begin{equation}\label{sec:Pa.10}
z^{*}_\zs{0}=\lambda_\zs{1}\k_\zs{*}(\Psi)
\quad\mbox{and}\quad
\tau^{*}=\lambda_\zs{2}\k_\zs{*}(\Psi)\,,
\end{equation}
where
$$
\lambda_\zs{1}=\max\left(
2c^{*}_\zs{1}M_\zs{1}\,,\,
2c^{*}_\zs{2}\,,\,
M_\zs{1}q^{*}\,,\,1\right)
\quad\mbox{and}\quad
\lambda_\zs{2}=\max\left(
c^{*}_\zs{1}M_\zs{1}\,,\,
c^{*}_\zs{2}
\right)\,.
$$
\noindent
Finally, we set
\begin{equation}\label{sec:Pa.11}
\gamma^{*}=\frac{1}{4\tau^{*}}\,.
\end{equation}

\medskip

\section{Continuous observations}\label{sec:Ci}

In this section we study the deviation in the ergodic theorem for the continuous observation 
case, which in this case is defined as
\begin{equation}\label{sec:Ci.1}
\Delta_\zs{T}(\phi)=\frac{1}{\sqrt{T}}\,
\int^T_0\,(\phi(y_\zs{t})\,-\,\pi_\zs{\vartheta}(\phi))\,\d t\,,
\end{equation}
where $\phi$ is any integrated function, i.e.   $|\phi|_\zs{1}<\infty$.

\begin{proposition}\label{Pr.sec:Ci.0}
For any $\nu_\zs{0}>0$ and $\nu_\zs{1}>0$

\begin{equation}\label{sec:Ci.1-00}
\sup_\zs{z\ge 0}
\quad e^{\kappa_\zs{0} z^{2}}
\quad
\sup_\zs{T\ge 1}
\quad
\sup_\zs{\phi\in \cV_\zs{\nu_\zs{0},\nu_\zs{1}}}
\quad
\sup_\zs{\vartheta\in\Theta}
\quad
\,
\P_\zs{\vartheta}\left(
|\Delta_\zs{T}(\phi)|\,
\ge z 
\right)\le 2
\,,
\end{equation}
where the parameter $\kappa_\zs{0}$ is given in \eqref{sec:Pa.5}.
\end{proposition}
\proof
Similarly to \cite{GaPe2} firstly we show that the deviation 
\eqref{sec:Ci.1} has an exponential moment, i.e. we show that
for the parameter $\kappa_\zs{0}$
\begin{equation}\label{sec:Ci.1-000}
\sup_\zs{T\ge 1}
\,
\sup_\zs{\vartheta\in\Theta}
\,
\E_\zs{\vartheta}
e^{\kappa_\zs{0}\Delta^{2}_\zs{T}(\phi)}\,
\le 2\,.
\end{equation}

\noindent Indeed, to show this inequality we need to estimate the expectation
of any even power for the deviation $\Delta_\zs{T}(\phi)$. To this end
we have to represent this deviation as the sum of a continuous
martingale and a negligible term. For this one needs to find 
a bounded solution for the following differential equation
\begin{equation}\label{sec:Ci.3}
\dot{v}_\zs{\vartheta}(u)\,+\,2\frac{S(u)}{\sigma^{2}(u)}v_\zs{\vartheta}(u)\, =\,2\frac{\wt{\phi}(u)}{\sigma^{2}(u)}, \quad
\wt{\phi}(u)=\phi(u)-\pi_\zs{\vartheta}(\phi)\,.
\end{equation}
One can check directly that the function
\begin{equation}\label{sec:Pr.4}
v_\zs{\vartheta}(u)=-2
\int^{\infty}_\zs{u}\,\frac{\wt{\phi}(y)}{\sigma^{2}(y)}\,
\exp\{2\int^{y}_\zs{u}\,S_\zs{1}(z)\d z\}\,\d y
\end{equation}
yields such a solution. We recall that the function
$S_\zs{1}$ is defined in \eqref{sec:Mr.2}. Moreover, due to Lemma~\ref{Le.sec:A.1}
from Appendix implies this function is uniform bounded.
By applying the Ito formula to the function
$V(y)=\int_0^y\,v_\zs{\vartheta}(u)\d u$ we following representation
\begin{equation}\label{sec:Pr.5}
\int^{T}_\zs{0}\wt{\phi}(y_\zs{s})\d s=V(y_\zs{T})
-V(y_\zs{0})
-\zeta_\zs{T}\,,
\end{equation}
where
$\zeta_\zs{T}=\int^{T}_\zs{0}\,v_\zs{\vartheta}(y_\zs{s})\sigma(y_\zs{s})\d w_\zs{s}$.
\noindent 
Therefore,  for any
 $T\ge 1$ through Lemma~\ref{Le.sec:A.1} we 
can estimate $\Delta_\zs{T}(\phi)$ from above
 as
$$
|\Delta_\zs{T}(\phi)|\,\le\,r|y_\zs{T}|\,+\,r|y_\zs{0}|\,+\,\frac{1}{\sqrt{T}}
\left|\zeta_\zs{T}\right|\,.
$$
\noindent 
Moreover, taking into account (see \cite{LpSh}, Lemma 4.11), that for any $m\ge 1$,
$$
\E_\zs{\vartheta}\,
\left(\zeta_\zs{T}\right)^{2m}
\le (2m-1)!!\,r^{2m}\sigma_\zs{\max}^{2m}\,T^m
\,,
$$ 
we obtain by Proposition~\ref{Pr.sec:A.3} , that for any $m\ge 1$
\begin{align*}
\E_\zs{\vartheta}|\Delta_\zs{T}(\phi)|^{2m}
&\le 3^{2m-1}
\left(r^{2m}(\E_\zs{\vartheta}|y_\zs{T}|^{2m}+|y_\zs{0}|^{2m})+
\,
\frac{\E_\zs{\vartheta}\,\left(\zeta_\zs{T}\right)^{2m}}{T^{m}}
\right)
\\[2mm]
&\le (3\,r)^{2m}
\left(4(m+1)(2m-1)!!\,\rho^{2m}+y^{2m}_\zs{0}+
(2m-1)!!\,\sigma_\zs{\max}^{2m}
\right)\,.
\end{align*}
Therefore,
 taking into account the definition of $\kappa_\zs{0}$, 
we obtain
\begin{align*}
\E_\zs{\vartheta}e^{\kappa_\zs{0}\Delta^{2}_\zs{T}(\phi)}&
=1\,+\,\sum_\zs{m=1}^{\infty}\frac{\kappa_\zs{0}^m}{m!}(3\,r)^{2m}
\left(4(2m+1)!!
\rho^{2m}+y^{2m}_\zs{0}+(2m-1)!!\sigma_\zs{\max}^{2m}\right)
\\[2mm]
&\le 1+\sum_\zs{m=1}^{\infty}\kappa_\zs{0}^{m}
(3\,r)^{2m}
\left(4(3\rho^{2})^{m}+y_\zs{0}^{2m}+2^m\sigma_\zs{\max}^{2m}\right)
\\[2mm]
&\le
1\,+\,\sum_\zs{m=1}^{\infty}(1/2)^{m}=2\,.
\end{align*}
From here we obtain the inequality \eqref{sec:Ci.1-000} and 
by the Chebychev inequality we come to the upper bound \eqref{sec:Ci.1-00}.
Hence Proposition~\ref{Pr.sec:Ci.0}.
 \endproof

\begin{remark}\label{Re.sec:Ci.1}
It should be noted that the inequality \eqref{sec:Ci.1-00}
is shown in \cite{GaPe2} for the process \eqref{sec:In.1} with 
$\sigma=1$. Thus Proposition~\ref{Pr.sec:Ci.0} extends teh result from 
\cite{GaPe2} for any diffusion function $\sigma$.
\end{remark}

\medskip

\section{Uniform geometric ergodicity}\label{sec:Ge}

Here we announce a result on geometric ergodicity obtained 
in \cite{GaPe4}.
\begin{theorem}\label{Th.sec:Ci.0}
There exist some constants $R\ge 1$ and 
$\kappa>0$ such that
\begin{equation}\label{sec:Ci.0}
\sup_\zs{t\ge 0} 
e^{\kappa t}
\sup_\zs{\|g\|_\zs{*}\le 1}
\sup_\zs{x\in\bbr}
\sup_\zs{\vartheta\in \Theta}\,
\frac{
\left|\E_\zs{\vartheta}\,\left( g(y_\zs{t})|y_\zs{0}=x\right)
-\pi_\zs{\vartheta}(g)\right|}{1+|x|}
\,\le R\,,
\end{equation} 
where the parameters $R$ and $\kappa$ are given in \cite{GaPe4}.
\end{theorem}

\medskip

\section{Burkh\"older's inequality}\label{Bi}

In this section we give the following inequality from 
\cite{DeDo},\cite{Ri}.
\begin{proposition}\label{Pr.sec:Bi.1}
 Let $(\Omega,\cF,(\cF_\zs{j})_\zs{1\le j\le n},\P)$
be a filtered probability space and
$(X_\zs{j}, \cF_\zs{j})_\zs{1\le j\le n}$ be  sequence of random
variables such that for some $p\ge 2$
$$
\max_\zs{1\le j\le n}\,\E\,|X_\zs{j}|^{p}\,<\,\infty\,.
$$
Define
$$
b_\zs{j,n}(p)=
\left(
\E\,
(
|X_\zs{j}|\,\sum^{n}_\zs{k=j}
|\E\,(X_\zs{k}|\cF_\zs{j})|
)^{p/2}
\right)^{2/p}\,.
$$
 Then
\begin{equation}\label{sec:Bi.1}
\E\,|
\sum^{n}_\zs{j=1}\,X_\zs{j}
\,|^{p}
\le\,(2p)^{p/2}
\left(
\sum^{n}_\zs{j=1}\,
b_\zs{j,n}(p)
\right)^{p/2}\,.
\end{equation}
\end{proposition}

\noindent
Proof of this Proposition is given in Appendix.

\medskip
 
\section{Proofs}\label{sec:Pr}

\subsection{Proof of Theorem~\ref{Th.sec:Mr.1}}

First note, that by Proposition~\ref{Pr.sec:A.3} and 
 the H\"older inequality we obtain for any $\alpha\ge 1$
\begin{equation}\label{sec:Pr.1}
\sup_\zs{t\ge 0}\,
\sup_\zs{\vartheta\in\Theta}\,
\E_\zs{\vartheta}\left(|y_\zs{t}|^{\alpha}|y_\zs{0}=x\right)\,\le\,
4\,(\alpha+1)^{\alpha/2}\,\rho^{\alpha}\,.
\end{equation}
Now we represent the deviation $D_\zs{T}(\phi)$ as
\begin{align}\nonumber
D_\zs{T}(\phi)
&=
\int^{T}_\zs{0}
(\phi(y_\zs{t})
-
\pi_\zs{\vartheta}(\phi))\d t 
+
\A_\zs{1,T}\,-\,\A_\zs{2,T}\\[2mm]\label{sec:Pr.1-1}
&=\sqrt{T}\,\Delta_\zs{T}(\phi)
+
\A_\zs{1,T}\,-\,\A_\zs{2,T}
\,,
\end{align}
where
$$
\A_\zs{1,T}=
\sum^{N}_\zs{j=1}
\int^{t_\zs{j}}_\zs{t_\zs{j-1}}\,
\left(
\phi(y_\zs{t_\zs{j}})
-
\phi(y_\zs{t})
\right)\,\d t
\quad\mbox{and}\quad
\A_\zs{2,T}=\int^{T}_\zs{\delta N}\,
(\phi(y_\zs{t})
-
\pi_\zs{\vartheta}(\phi))\d t
\,. 
$$
To estimate the term $\A_\zs{1,T}$ 
we represent through the Ito formula the difference
$\phi(y_\zs{t_\zs{j}})
-
\phi(y_\zs{t})$ as 
\begin{align*}
\phi(y_\zs{t_\zs{j}})
-
\phi(y_\zs{t})
&=
\int^{t_\zs{j}}_\zs{t}\,
\cL_\zs{\vartheta}(\phi)(y_\zs{s})\,\d s
+
\int^{t_\zs{j}}_\zs{t}
\dot{\phi}(y_\zs{s})\sigma(y_\zs{s})\d W_\zs{s}\\
&=
\wt{\pi}_\zs{\vartheta}(\phi)
(t_\zs{j}-t)
+
\Psi_\zs{j}(t)
+
\int^{t_\zs{j}}_\zs{t}
\dot{\phi}(y_\zs{s})\sigma(y_\zs{s})\d W_\zs{s}
\,,
\end{align*}
where 
$$
\Psi_\zs{j}(t)
=
\int^{t_\zs{j}}_\zs{t}\,
\psi(y_\zs{s})\,\d s\,,
\quad
\omega_\zs{j}(t)=
\int^{t_\zs{j}}_\zs{t}\,
\dot{\phi}(y_\zs{s})\sigma(y_\zs{s})\,\d W_\zs{s}
$$
and $\psi(y)=\cL_\zs{\vartheta}(\phi)(y)-\wt{\pi}_\zs{\vartheta}(\phi) $.
Now setting
$$
X_\zs{j}=
\int^{t_\zs{j}}_\zs{t_\zs{j-1}}
\Psi_\zs{j}(t)
\,\d t
\quad\mbox{and}\quad
\eta_\zs{j}=
\int^{t_\zs{j}}_\zs{t_\zs{j-1}}
\omega_\zs{j}(t)
\d t\,,
$$
we obtain
\begin{equation}\label{sec:Pr.2}
\A_\zs{1,T}=\wt{\pi}_\zs{\vartheta}(\phi)\frac{N\delta^{2}}{2}+
\sum^{N}_\zs{j=1}\,X_\zs{j}
+
\sum^{N}_\zs{j=1}\,\eta_\zs{j}\,.
\end{equation}

To estimate the second term in the right-hand part of \eqref{sec:Pr.2}, 
we make use of the Proposition~\ref{Pr.sec:Bi.1}.We start with verifying its conditions. 
Putting $\cF_\zs{s}=\sigma\{y_\zs{u}\,,0\le u\le s\}$, 
we obtain by Theorem~\ref{Th.sec:Ci.0}, that  for any $t\ge s$
and for any $\phi$ from the functional class \eqref{sec:Mr.2-2}
$$
|\E_\zs{\vartheta} \left(
\psi(y_\zs{t})
|\cF_\zs{s}
\right)|
\le \,\mu(\phi) \,
R\, 
\left(
1+|y_\zs{s}|
\right)\,
e^{-\kappa (t-s)}
\le \,\nu_\zs{3}\,
R\, 
\left(
1+|y_\zs{s}|
\right)\,
e^{-\kappa (t-s)}
\,.
$$

\noindent 
Therefore, for any  $k> j$,
\begin{equation}\label{sec:Pr.3}
|\E_\zs{\vartheta} (X_\zs{k}|\cF_\zs{t_\zs{j}})|
\le 
R e^{\kappa}
(
1+
|y_\zs{t_\zs{j}}|
)
\,
\nu_\zs{3}
\,\delta^{2}
\,e^{-\kappa \delta(k-j)}
\,.
\end{equation}
It should be noted also, that the random variables $X_\zs{j}$ are bounded, i.e.
$$
|X_\zs{j}|\le \nu_\zs{3}\delta^{2}\,.
$$
To estimatie the probability tail for the sum 
$\sum^{n}_\zs{j=1} X_\zs{j}$ we will use the inequality \eqref{sec:Bi.1}.
For this we need to estmate the coefficients
$b_\zs{j,N}(p)$ for any 
$p\ge 1$. From here, taking into account that 
$1-e^{-\kappa\delta}\ge \kappa\delta e^{-\kappa}$
and that for $p\ge 2$
$$
\left(
\E_\zs{\vartheta}(1+
|y_\zs{t_\zs{j}}|)^{p/2}
\right)^{2/p}
\le
1+
\left(\E_\zs{\vartheta}|y_\zs{t_\zs{j}}|^{p/2}\right)^{2/p}\,,
$$
we can estimate the coefficient $b_\zs{j,N}(p)$ as
$$
b_\zs{j,N}(p)\,\le\,\frac{1}{\kappa}\,R\,e^{2\kappa}\,
\varsigma^{2}
\,\left(1+
(\E|y_\zs{t_\zs{j}}|^{p/2})^{2/p}
\right)\,,
$$
where $\varsigma^{2}=\nu^{2}_\zs{3}\delta^{3}$. Now the inequality \eqref{sec:Pr.1} yields
$$
b_\zs{j,N}(p)\,
\le 
R_\zs{1}
\varsigma^{2}
\sqrt{2+p}
\le R_\zs{1}
\varsigma^{2}
\sqrt{2p}
\,,
$$
where 
$$
R_\zs{1}= \frac{1}{\kappa}\,R\,e^{2\kappa}(1+\rho)
\,.
$$
Using this in \eqref{sec:Bi.1}
 we obtain, that 
for any $p>2$,
\begin{align*}
\E_\zs{\vartheta}\,|\sum^{N}_\zs{k=1}\,X_\zs{k}|^{p}
&\le (2p)^{p/2}\, N^{p/2}\,
R_\zs{1}^{p/2}\,\varsigma^{p}\,(2p)^{p/4}\\[2mm]
&\le (2\sqrt{R_\zs{1}}\,\varsigma)^{p}\,N^{p/2}\,p^{p}\,.
\end{align*}
Therefore, by Chebyshev's inequality 
$$
\P_\zs{\vartheta}\left(
|\sum^{N}_\zs{k=1}\,X_\zs{k}|
\ge z\sqrt{N}
\right)
\le 
e^{p\ln(a)+ p\ln p}
$$
with $a=2\sqrt{R_\zs{1}}\varsigma/z$. 
Minimizing now the right-hand part over $p\ge 2$,
 we obtain for  $z\ge 4e\sqrt{R_\zs{1}}\nu_\zs{3}\delta^{3/2}$
\begin{equation}\label{sec:Pr.4-0}
\P_\zs{\vartheta}\left(
|\sum^{N}_\zs{k=1}\,X_\zs{k}|
\ge z\sqrt{N}
\right)
\le\,e^{-z/\varsigma_\zs{1}}\,,
\end{equation}
 where $\varsigma_\zs{1}=2e\sqrt{R_\zs{1}}\varsigma$.

Moreover, note that by the Burkholder-Davis-Gundy inequality, for any $\alpha\ge 1$,
$$
\E_\zs{\vartheta}\,
|\omega_\zs{j}(t)|^{\alpha}
\le\,
(\alpha)^{\alpha/2}\,
\nu^{\alpha}_\zs{2}
\sigma_\zs{\max}^{\alpha}\,(t_\zs{j}-t)^{\alpha/2}
\,.
$$
Using this and the the H\"older inequality,
we get 
$$
\E_\zs{\vartheta}\,|\eta_\zs{j}|^{\alpha}\,
\le
\delta^{\alpha-1}\,\int^{t_\zs{j}}_\zs{t_\zs{j-1}}\,
\E_\zs{\vartheta}\,
|\omega_\zs{j}(t)|^{\alpha}\,\d t
\le\,
\,\delta^{3\alpha/2}\,\alpha^{\alpha/2}\,
\nu^{\alpha}_\zs{2}
\sigma_\zs{\max}^{\alpha}\,.
$$
Note, that in this case in the right hand of the inequality
\eqref{sec:Bi.1} 
$$
b_\zs{j,N}=\left(
\E_\zs{\vartheta}\,|\eta_\zs{j}|^{p}
\right)^{2/p}\,.
$$
Therefore, similarly to the inequality \eqref{sec:Pr.4}
 we find, that for all $z\ge 2\varsigma_\zs{2}$,
\begin{equation}\label{sec:Pr.5-0}
\P_\zs{\vartheta}\left(
|\sum^{N}_\zs{k=1}\,\eta_\zs{k}|
\ge z\sqrt{N}
\right)
\le 
\,e^{-z/\varsigma_\zs{2}}\,,
\end{equation}
where 
$\varsigma_\zs{2}=\sqrt{2}e\delta^{3/2}\nu_\zs{2}
\sigma_\zs{\max}$.
\noindent 
Now from \eqref{sec:Pr.2}, \eqref{sec:Pr.4}--\eqref{sec:Pr.5} it follows that
for $z\ge z_\zs{0}$
\begin{align}\nonumber
\P_\zs{\vartheta}\left(|\A_\zs{1,T}|\ge z\sqrt{N}\right)
&\le
\,
\P_\zs{\vartheta}\left(
|\sum^{N}_\zs{k=1}\,X_\zs{k}|
\ge z\sqrt{N}/4\right)\\[2mm]\label{sec:Pr.6}
&+
\P_\zs{\vartheta}\left(
|\sum^{N}_\zs{k=1}\,\eta_\zs{k}|
\ge z\sqrt{N}/4\right)
\le
2\,e^{-z/4\tau}
\,,
\end{align}
when the parameters 
$z_\zs{0}$ and $\tau$ are given in \eqref{sec:Pa.6}. 
Moreover, note that due to \eqref{sec:Mr.2-2} the
last term 
in \eqref{sec:Pr.1-1} is bounded, i.e.
$$
|\A_\zs{2,T}|\le 2\delta \|\phi\|_\zs{*}
\le 2\delta \nu_\zs{2}
\le z_\zs{0}\sqrt{N}/4\,. 
$$
Finally, from  \eqref{sec:Pr.1-1} for $z\ge z_\zs{0}$
 one has
\begin{align*}
\P_\zs{\vartheta}(|D_\zs{T}(\phi)|\ge z\sqrt{N})&\le
\P_\zs{\vartheta}\left(\sqrt{T}|\Delta_\zs{T}(\phi)|
\,+\,
|\A_\zs{1,T}|\ge 3z\sqrt{N}/4\right)
\\[2mm]
&\le
\P_\zs{\vartheta}\left(\sqrt{T}|\Delta_\zs{T}(\phi)|\ge 3z\sqrt{N}/8\right)
\\[2mm]
&+
\P_\zs{\vartheta}\left(|\A_\zs{1,T}|\ge 3z\sqrt{N}/8\right)
\,. 
\end{align*}

\noindent 
Taking into account here, that $N/T\ge (1-\delta)/\delta$
for any $0<\delta<1$ and $T\ge 1$, we obtain, that
\begin{align*}
\P_\zs{\vartheta}(|D_\zs{T}(\phi)|\ge z\sqrt{N})
&\le
\P_\zs{\vartheta}\left(|\Delta_\zs{T}(\phi)|\ge 
\frac{3z\sqrt{(1-\delta)}}{8\sqrt{\delta}}\right)
\\[2mm]
&+\,
\P_\zs{\vartheta}\left(|\A_\zs{1,T}|\ge \frac{3}{8}z\sqrt{N}\right)
\,. 
\end{align*}
Therefore, applying here the inequalities
\eqref{sec:Ci.1-00} and \eqref{sec:Pr.6}
we come to the upper bound \eqref{sec:Mr.3}
with the parameter $\varkappa$ given in \eqref{sec:Pa.7}.
Hence Theorem~\ref{Th.sec:Mr.1}.
\endproof

\medskip

\subsection{Proof of Theorem~\ref{Th.sec:Mr.2}}

Firstly, note that in this case
$$
|\psi_\zs{h,x_\zs{0}}|_\zs{1}=|\Psi|_\zs{1}\,,\quad
\|\psi_\zs{h,x_\zs{0}}\|_\zs{*}=\frac{1}{h}\|\Psi\|_\zs{*}
\quad\mbox{and}\quad
\|\dot{\psi}_\zs{h,x_\zs{0}}\|_\zs{*}=\frac{1}{h^{2}}\|\dot{\Psi}\|_\zs{*}\,.
$$

\noindent
Moreover, taking into account  that $|S(y)|\le M+L\x_\zs{*}+L|y|$, we find
that
\begin{equation}\label{sec:Pr.8}
\sup_\zs{|y|\le |x_\zs{0}|+2}|S(y)|
\le 
M_\zs{1}\,,
\end{equation}
where
$M_\zs{1}$ is given in \eqref{sec:Pa.8}.

\noindent
Therefore, in view of the fact that $0<h<1$, we can estimate from above 
the parametrs \eqref{sec:Pa.1}
as
\begin{equation}\label{sec:Pr.9}
\mu(\psi_\zs{h,x_\zs{0}})\le \mu_\zs{*} h^{-3}
\quad\mbox{and}\quad
\wt{\mu}(\psi_\zs{h,x_\zs{0}})\le \wt{\mu}_\zs{*} h^{-2}
\,,
\end{equation}
where
$$
\mu_\zs{*}= 
\,\max\left(
\|\dot{\Psi}\|_\zs{*}
\,,\,
\|\ddot{\Psi}\|_\zs{*}
\right)
\,M_\zs{1}
\quad\mbox{and}\quad
 \wt{\mu}_\zs{*}= 
\,\max\left(
|\dot{\Psi}|_\zs{1}
\,,\,
|\ddot{\Psi}|_\zs{1}
\right)
\,M_\zs{1}q^{*}
\,.
$$
\noindent
Therefore, the function $\psi_\zs{h,x_\zs{0}}$ belongs to the 
class \eqref{sec:Mr.2-2} with the following parameters
$$
\nu_\zs{0}=|\Psi|_\zs{1}\,,\quad
\nu_\zs{1}=\frac{\|\Psi\|_\zs{*}}{h}\,,\quad
\nu_\zs{2}=\frac{\|\dot{\Psi}\|_\zs{*}}{h^{2}}\,,\quad
\nu_\zs{3}=\frac{\mu_\zs{*}}{h^{3}}\,,\quad
\nu_\zs{4}=\frac{\wt{\mu}_\zs{*}}{h^{2}}\,.
$$
Therefore, in this case the coefficient \eqref{sec:Pa.5}
equals to $\kappa_\zs{0}(|\Psi|_\zs{1})$ and
the parameters \eqref{sec:Pa.6} can be represented as
\begin{align}\nonumber
z_\zs{0}&=\frac{\delta^{3/2}}{h^{3}}\,
\max\left(
2c^{*}_\zs{1}\mu_\zs{*}\,,\,
2c^{*}_\zs{2}\|\dot{\Psi}\|_\zs{*}h\,,\,
\wt{\mu}_\zs{*}hT^{1/2}\,,\,
\|\Psi\|_\zs{*}h^{2}T^{-1/2}
\right)
\\[2mm]\label{sec:Pr.10}
\tau&=\frac{\delta^{3/2}}{h^{3}}\,
\max\left(
c^{*}_\zs{1}\mu_\zs{*}\,,\,
c^{*}_\zs{2}\|\dot{\Psi}\|_\zs{*}h
\right)\,.
\end{align}
Therefore, thanks to the condition \eqref{sec:Mr.5}  for any 
$T^{-1/2}\le h\le 1$
\begin{equation}\label{sec:Pr.7}
 z_\zs{0}\le l^{-3/2}_\zs{T}\,z^{*}_\zs{0}
\quad\mbox{and}\quad
\tau
\le 
l^{-3/2}_\zs{T}\,\tau^{*}\,,
\end{equation}
where the parameters $z^{*}_\zs{0}$ 
and $\tau^{*}$ are given in \eqref{sec:Pa.10}.
Note now that, by the condition \eqref{sec:Mr.4}
$$
\P_\zs{\vartheta}\left(
|D_\zs{T}(\psi_\zs{h,x_\zs{0}})|\,
\ge
\,a\,T
\right)
\le 
\P_\zs{\vartheta}\left(
|D_\zs{T}(\psi_\zs{h,x_\zs{0}})|\,
\ge
\,z_\zs{1}\,\sqrt{N}
\right)
$$
where $z_\zs{1}=a/\sqrt{l_\zs{T}}$. The first inequality in
\eqref{sec:Pr.7} implies that 
$z_\zs{1}\ge z_\zs{0}$ for all $a\ge a_\zs{*}=z^{*}_\zs{0}/l_\zs{T}$.
Moreover, from the last inequality in 
\eqref{sec:Pr.7} it follows, that for $a\ge a_\zs{*}$
$$
\min\left(\varkappa z_\zs{1}\,,\,\gamma\right)=
\min\left(
\varkappa z_\zs{1}\,,\,\frac{1}{4\tau}
\right)
\ge
\min\left(\varkappa \frac{z^{*}_\zs{0}}{l_\zs{T}\sqrt{l_\zs{T}}}\,,\,
\frac{l_\zs{T}\sqrt{l_\zs{T}}}{4\tau^{*}}\right)\,.
$$
Taking into account here the definition of $\varkappa$
in \eqref{sec:Pa.7} and the  form for $\delta$ given by 
\eqref{sec:Mr.4} we obtain that for sufficiently large $T$
$$
\min\left(\varkappa \frac{z^{*}_\zs{0}}{l_\zs{T}\sqrt{l_\zs{T}}}\,,\,
\frac{l_\zs{T}\sqrt{l_\zs{T}}}{4\tau^{*}}\right)=
\frac{l_\zs{T}\sqrt{l_\zs{T}}}{4\tau^{*}}\,.
$$
Thus, through Theorem~\ref{Th.sec:Mr.1} 
 we come to the inequality
 \eqref{sec:Mr.7}. 
Hence Theorem~\ref{Th.sec:Mr.2}
 \endproof

\subsection{Proof of Theorem~\ref{Th.sec:Mr.3}}

First we represent the tail probability as
$$
\P_\zs{\vartheta}\left(
|D_\zs{T}(\chi_\zs{h,x_\zs{0}})|\,
\ge\,
\epsilon_\zs{T}\,T
\right)
=
\I_\zs{1}
+
\I_\zs{2}\,,
$$
where
$$
\I_\zs{1}=\P_\zs{\vartheta}\left(
\sum^{N}_\zs{j=1}\chi_\zs{h,x_\zs{0}}(y_\zs{t_\zs{j}})\Delta t_\zs{j}
\le
(\pi_\zs{\vartheta}(\chi_\zs{h,x_\zs{0}})-\epsilon_\zs{T})\,T
\right)
$$
and
$$
\I_\zs{2}=\P_\zs{\vartheta}\left(
\sum^{N}_\zs{j=1}\chi_\zs{h,x_\zs{0}}(y_\zs{t_\zs{j}})\Delta t_\zs{j}
\ge
(\pi_\zs{\vartheta}(\chi_\zs{h,x_\zs{0}})+\epsilon_\zs{T})\,T
\right)\,.
$$
Let us define now the following smoothing indicator functions
$$
\Psi_\zs{1,\eta}(u)=\frac{1}{\eta}
\int^{+\infty}_\zs{-\infty}\,
\Chi_\zs{\{|z|\le 1-\eta\}}\,V\left(\frac{z-u}{\eta}\right)
\d z
$$
and
$$
\Psi_\zs{2,\eta}(u)=\frac{1}{\eta}
\int^{+\infty}_\zs{-\infty}\,
\Chi_\zs{\{|z|\le 1+\eta\}}\,V\left(\frac{z-u}{\eta}\right)
\d z
\,,
$$
where $\eta$ is a smoothing positive parameter which will be specified later,
$V$ 
 is a two times continuously differentiable even
$\bbr\to\bbr$ function such that $V(z)=0$ for $|z|\ge 1$ and
$$
\int^{1}_\zs{-1}\,V(z)\d z=1\,.
$$
It is easy to see that, for any $y\in\bbr$ and $0<\eta\le 1/2$,
$$
\Psi_\zs{1,\eta}(u)(y)\le
\chi(y)
\le \Psi_\zs{2,\eta}(y)
$$
and $\Psi_\zs{2,\eta}(y)=0$ for $|y|\ge 2$.
Moreover, for the functions
$$
\psi_\zs{i,h}(y)=\frac{1}{h}\,\Psi_\zs{i,\eta}
\left(\frac{y-z_\zs{0}}{h}\right)
$$
using the inequality \eqref{sec:A.5},
we can estimate the difference between the cooresponding
ergodic intergals \eqref{sec:In.4}  as
$$
|\pi_\zs{\vartheta}(\chi_\zs{h,x_\zs{0}})
-
\pi_\zs{\vartheta}(\psi_\zs{i,h})
|
\le 4\eta q^{*}\,.
$$
\noindent 
Therefore, choosing here $\eta=\epsilon^{2}_\zs{T}$
we obtain, for sufficiently large $T$,
$$
\I_\zs{i}\le 
\P_\zs{\vartheta}\left(
|D_\zs{T}(\phi_\zs{i,h})|
\ge
\,\epsilon_\zs{T}\,T/2
\right)\,.
$$
One can check directly that in this case the operator
\eqref{sec:Pa.9} has the following asymptotic ($T\to\infty$) form 
$$
\k_\zs{*}(\Psi_\zs{i,\eta})=\aO\left(\eta^{-2}\right)\,.
$$

\noindent Therefore, from 
\eqref{sec:Pa.10} and \eqref{sec:Pr.7}
it follows that
for $T\to\infty$ and $h\ge T^{-1/2}$  
$$
z_\zs{0}(\phi_\zs{i,h})=
\aO\left(\eta^{-2} l^{-3/2}_\zs{T}\right)
\quad\mbox{and}\quad
\tau(\phi_\zs{i,h})=\aO\left(\eta^{-2} l^{-3/2}_\zs{T}\right)
\,,
$$
i.e.
$$
z_\zs{0}(\phi_\zs{i,h})=\aO\left(\frac{1}{\epsilon^{4}_\zs{T} l^{3/2}_\zs{T}}\right)
\quad\mbox{and}\quad
\tau(\phi_\zs{i,h})=\aO\left(\frac{1}{\epsilon^{4}_\zs{T} l^{3/2}_\zs{T}}\right)\,.
$$
Now we have
$$
\P_\zs{\vartheta}\left(
|D_\zs{T}(\psi_\zs{h,x_\zs{0}})|\,
\ge
\,\epsilon_\zs{T} T
\right)
\le 
\P_\zs{\vartheta}\left(
|D_\zs{T}(\psi_\zs{h,x_\zs{0}})|\,
\ge
\,z_\zs{1}\,\sqrt{N}
\right)\,,
$$
where $z_\zs{1}=\epsilon_\zs{T}/\sqrt{l_\zs{T}}$.
The last equality in \eqref{sec:Mr.6} implies 
$z_\zs{1}\ge z_\zs{0}$ for sufficiently large $T$.
Moreover, taking into account, that
there exists a constant $\c_\zs{*}>0$ such that
 for sufficiently large $T$
$$
\varkappa z_\zs{1} \ge \c_\zs{*} T\sqrt{l_\zs{T}}\epsilon_\zs{T}
\quad\mbox{and}\quad
\gamma\ge \c_\zs{*} l_\zs{T}\sqrt{l_\zs{T}}\epsilon^{4}_\zs{T}\,,
$$
i.e. for sufficiently large $T$
$$
\min\left(
\varkappa z_\zs{1} \,,\,
\gamma
\right)
\ge c_\zs{*}
l_\zs{T}\sqrt{l_\zs{T}}\epsilon^{4}_\zs{T}\,.
$$
Therefore, by Theorem~\ref{Th.sec:Mr.1} for sufficiently large $T$
$$
\P_\zs{\vartheta}\left(
|D_\zs{T}(\psi_\zs{h,x_\zs{0}})|\,
\ge
\,\epsilon_\zs{T} T
\right)
\le 4
e^{-\c_\zs{*} l_\zs{T}\epsilon^{5}_\zs{T}}
\,.
$$

\noindent 
Now the last condition in \eqref{sec:Mr.6}
yields the equality \eqref{sec:Mr.8}. Hence 
Theorem~\ref{Th.sec:Mr.3}.
\endproof

\medskip

\setcounter{section}{0}
\renewcommand{\thesection}{\Alph{section}}

\section{Appendix}\label{sec:A}

\subsection{Proof of Proposition~\ref{Pr.sec:Bi.1}}

We set
$$
h_\zs{n}(t)=\E |S_\zs{n-1}+tX_\zs{n}|^{p}
\quad\mbox{with}\quad S_\zs{n}=\sum^n_\zs{j=1}X_\zs{j}\,.
$$
By the induction method we assume that for any $1\le k\le n-1$ and $0\le t\le 1$
\begin{equation}\label{sec:A.2}
h_\zs{k}(t)\le 
(2p)^{p/2}\,
B_\zs{k}^{p/2}(t)\,,
\end{equation}
where
$$
B_\zs{k}(t)=\sum^{k-1}_\zs{j=1}\,b_\zs{j,k}(p)
+
t
b_\zs{k,k}(p)\,.
$$
Note now that as is shown in \cite{Ri} (Theorem 2.3)
\begin{equation}\label{sec:A.3}
\E |S_\zs{n}|^{p}=
p(p-1)
\sum^n_\zs{j=1}
\int^{1}_\zs{0}\,\E |S_\zs{j-1}+vX_\zs{j}|^{p-2}
(
-vX2_\zs{j}+\Upsilon(j,n)
)\d v\,.
\end{equation}
with
$$
\Upsilon(j,n)=
X_\zs{j}\sum^n_\zs{k=j}\E(X_\zs{k}|\cF_\zs{j})\,.
$$
Therefore,
\begin{align*}
 h_\zs{n}(t)&=p(p-1)
\sum^{n-1}_\zs{j=1}
\int^{1}_\zs{0}\,\E |S_\zs{j-1}+vX_\zs{j}|^{p-2}
(
-vX2_\zs{j}+G(i,n,t)
)\d v\\[2mm]
&+
p(p-1)
\int^{1}_\zs{0}\,\E |S_\zs{n-1}+vtX_\zs{n}|^{p-2}
t^{2}(
1-v
)X^{2}_\zs{n}\d v
\,,
\end{align*}
where
$$
G(j,n,t)=
\Upsilon(j,n-1)
+
t
X_\zs{j}\,
\E (X_\zs{n}|\cF_\zs{j})
\,.
$$
Moreover, we can estimate $h_\zs{n}(t)$ as
\begin{align*}
\frac{h_\zs{n}(t)}{p2}&\le 
\sum^{n-1}_\zs{j=1}
\int^{1}_\zs{0}\,\E |S_\zs{j-1}+vX_\zs{j}|^{p-2}\,|G(i,n,t)|\d v\\[2mm]
&+
\int^{t}_\zs{0}\,\E |S_\zs{n-1}+sX_\zs{n}|^{p-2}
\,X^{2}_\zs{n}\d s
\end{align*}
Now taking into account that for $0\le t\le 1$
$$
\left(
\E |G(j,n,t)|^{p/2}
\right)^{2/p}
\,\le\,b_\zs{j,n}(p)\,,
$$
we obtain by the H\"older inequality
$$
\int^{1}_\zs{0}\,
\E |S_\zs{j-1}+vX_\zs{j}|^{p-2}
 \,|G(i,n,t)|\,\d v\,\le \,
\int^{1}_\zs{0}\,h^{\alpha}_\zs{j}(v)\, b_\zs{j,n}(p)\d v\,,
$$
where $\alpha=1-2/p$. 
Therefore,
$$
\frac{h_\zs{n}(t)}{p2}\le 
\sum^{n-1}_\zs{j=1}\,b_\zs{j,n}(p)
\int^{1}_\zs{0}\,h^{\alpha}_\zs{j}(v)\d v
+b_\zs{n,n}(p)
\int^{t}_\zs{0}\,h^{\alpha}_\zs{n}(s)\d s
$$

\noindent

Now by the induction assumption
for any $1\le j\le n-1$
$$
b_\zs{j,n}(p)
\int^{1}_\zs{0}\,h^{\alpha}_\zs{j}(v)\d v
\,
\le \,(2p)^{(p-2)/2}\,
\int^{1}_\zs{0}\,
B^{(p-2)/2}_\zs{j}(v)\,\d v\, b_\zs{j,n}(p)\,.
$$
Moreover, taking into account that 
$$
B_\zs{j}(v)\le \sum^{j-1}_\zs{i=1}\,b_\zs{i,n}+vb_\zs{j,n}(p)\,,
$$
we obtain that
$$
\int^{1}_\zs{0}\,
B^{(p-2)/2}_\zs{j}(v)\,\d v\, b_\zs{j,n}(p)\,\le\,
\frac{2}{p}
\left(
(\sum^{j}_\zs{i=1}\,b_\zs{i,n})^{p/2}
-
(\sum^{j-1}_\zs{i=1}\,b_\zs{i,n})^{p/2}
\right)\,.
$$
This implies for any $0\le t\le 1$
\begin{equation}\label{sec:A.4}
h_\zs{n}(t)\le k_\zs{n}\,\int^t_\zs{0}\,h^{\alpha}_\zs{n}(v)\,\d v
+ f_\zs{n}
\end{equation}
with
$$
k_\zs{n}=p^2b_\zs{n,n}(p)
\quad\mbox{and}\quad
f_\zs{n}=\left(2p\sum^{n-1}_\zs{j=1} b_\zs{j,n}(p)\right)^{p/2}\,.
$$
Now by setting
$$
Z(t)=\int^t_\zs{0}\,h^{\alpha}_\zs{n}(s)\d s
+\frac{f_\zs{n}}{k_\zs{n}}\,,
$$
we obtain from \eqref{sec:A.4} that
$$
\dot{Z}(t)\le k^{\alpha}_\zs{n} Z^{\alpha}(t)\,.
$$
Now introducing
$$
g(t)=\dot{Z}(t)- k^{\alpha}_\zs{n} Z^{\alpha}(t)\,,
$$
we obtain the differential equation
$$
\dot{Z}(t)= k^{\alpha}_\zs{n} Z^{\alpha}(t)
+
g(t)
$$
with $g(t)\le 0$.
From here we obtain 
$$
Z^{2/p}(t)=Z^{2/p}(0)+\frac{2}{p}k^{\alpha}_\zs{n}t+
\int^t_\zs{O}\frac{g(u)}{Z^{\alpha}(u)}\,\d u
\le\,
Z^{2/p}(0)+\frac{2}{p}k^{\alpha}_\zs{n}t
\,,
$$
i.e.
$$
Z(t)\le\,
\left(
Z^{2/p}(0)+\frac{2}{p}k^{\alpha}_\zs{n}t
\right)^{p/2}\,.
$$
Substituting this bound in \eqref{sec:A.4} we obtain
\begin{align*}
h_\zs{n}(t)&\le k_\zs{n}Z(t)\le 
k_\zs{n}
\left(
Z^{2/p}(0)+\frac{2}{p}k^{\alpha}_\zs{n}t
\right)^{p/2}\\[2mm]
&=
\left(
2p\sum^{n-1}_\zs{j=1} b_\zs{j,n}(p)
+
2 p t b_\zs{n,n}(p)
\right)^{p/2}\,.
\end{align*}
Hence Proposition~\ref{Pr.sec:Bi.1}. \endproof

\subsection{Uniform bound for the invariant density }

\begin{lemma}\label{Le.sec:A.2}  The invariant density 
\eqref{sec:Mr.2} is uniformly bounded: 
\begin{equation}\label{sec:A.5}
\sup_\zs{x\in\bbr}\sup_\zs{\vartheta\in\Theta}\,q_\zs{\vartheta}(x)\,
\le\,q^{*}<\infty\,,
\end{equation}
where the upper bound $q^{*}$ is given in \eqref{sec:Pa.3}.
\end{lemma}
\proof
First, note that through the definition of $\Theta$ 
we can check directly that for any $|x|\ge x_\zs{*}$
\begin{equation}\label{sec:A.6}
2\int^{x}_\zs{0}\,S_\zs{1}(v)\,d\,v
\le \beta_\zs{1}|x|-\beta_\zs{2}(|x|-\x_\zs{*})^{2}\,,
\end{equation}
where the coefficients $\beta_\zs{1}$ and $\beta_\zs{2}$ are given in 
\eqref{sec:Pa.2}. Therefore, taking into account, that for  $|x|\ge x_\zs{*}$
$$
2\int^{x}_\zs{0}\,S_\zs{1}(v)\,d\,v\le \beta_\zs{1}
\,x_\zs{*}\,,
$$
we obtain that
$$
2\sup_\zs{x\in\bbr}\,
\int^{x}_\zs{0}\,S_\zs{1}(v)\,d\,v
\le
\beta_\zs{1}x_\zs{*}+
\frac{\beta_\zs{1}}{4\beta_\zs{2}}\,.
$$
Estimating now the denominator in \eqref{sec:Mr.2}
from below as
$$
\int_\zs{\bbr}\sigma^{-2}(z)\,
e^{\wt{S}(z)}\d z\,
\ge \int^{1}_\zs{0}\sigma^{-2}(z)\,\d z\ge 
\frac{1}{\sigma^{2}_\zs{max}}\,,
$$
and taking into account the definition of $q^{*}$ 
we come to the upper uniform bound \eqref{sec:A.5}.
Hence Proposition~\ref{Le.sec:A.2}. 
\endproof
 
\medskip

\subsection{Moment bound for the process $y_\zs{t}$. }

\begin{proposition}\label{Pr.sec:A.3} For any $m\ge 1$
$$
\sup_\zs{t\ge 0}
\,\sup_\zs{\vartheta\in\Theta}\,\E_\zs{\vartheta}|y_\zs{t}|^{2m}\,\le\,
4(m+1)(2m-1)!!\,\rho^{2m}
\le 4(2m)^{m}\,\rho^{2m}\,,
$$
where $\rho$ is given in \eqref{sec:Pa.5}.
\end{proposition}

\proof
First note, that through the Ito formula
we can write for the function $z_\zs{t}(m)=\E_\zs{\vartheta}y_\zs{t}^{2m}$
the following intergal equality
\begin{align*}
z_\zs{t}(m)&=z_\zs{0}(m)+2m\int^{t}_\zs{0}\E_\zs{\vartheta}y_\zs{s}^{2m-1}\,
S(y_\zs{s})\d s\\[2mm]
&+
m(2m-1)\int^{t}_\zs{0}\E_\zs{\vartheta}y_\zs{s}^{2m-2}\sigma^{2}(y_\zs{s})ds\,,
\end{align*}
which can be rewritten as the differential equality
$$
\dot{z}_\zs{t}(m)=2m \E_\zs{\vartheta}y_\zs{s}^{2m-1}\,
S(y_\zs{t})
+
m(2m-1) \E_\zs{\vartheta}y_\zs{t}^{2m-2}\sigma^{2}(y_\zs{t})\,.
$$
\noindent 
Taking into account here that
$\sup_\zs{x\in\bbr}\sigma^{2}(x)\le \sigma^{2}_\zs{\max}$ 
we obtain, that for any $m\ge 1$ and $t\ge 0$
$$
\dot{z}_\zs{t}(m)\le 2m\E_\zs{\vartheta}y_\zs{t}^{2m-1}S(y_\zs{t})
+m(2m-1)\sigma_\zs{\max}^2z_\zs{t}(m-1)\,. 
$$
Now we need to estimate from above 
the function $x^{2m-1}S(x)$. Obviously, that for any $K>\x_\zs{*}$
$$
x^{2m-1}S(x)\le K^{2m-1}\sup_\zs{|x|\le K}|S(x)|
\Chi_\zs{\{|x|\le K\}}\,+\,x^{2m}\frac{S(x)}{x}\Chi_\zs{\{|x|>K\}}\,.
$$
 Taking into account that $\sup_\zs{|x|>\x_\zs{*}}|\dot{S}(x)|\le L$, we obtain, for any
$x\in [\x_\zs{*},K]$,
$$
|S(x)|\le |S(\x_\zs{*})|+L |x-\x_\zs{*}|
\le M+L(K-\x_\zs{*})\,.
$$
Similarly, we obtain the same upper bound for $x\in [-K,\,-\x_\zs{*}]$. Therefore,
$$
\sup_\zs{|x|\le K}|S(x)|\le M\,+\,L\,(K-\x_\zs{*}).
$$
Consider now the case $|x|>K$. 
We recall, that $\sup_\zs{|x|\ge\x_\zs{*}}\dot{S}(x)\le -L^{-1}$. Therefore, 
$$
\frac{S(x)}{x}\le \frac{M}{K}-
\frac{K-\x_\zs{*}}{L K}\,.
$$
Choosing   $K=2(\x_\zs{*}+ML)$ yields
$$
\frac{S(x)}{x}\le-\frac{1}{2L}\,.
$$
Therefore,
\begin{align*}
x^{2m-1}S(x)&\le 
K^{2m-1}\left(M+L(K-\x_\zs{*})\right)
-\frac{1}{2L}x^{2m}\Chi_\zs{\{|x|>K\}}
\\[2mm]
&
=K^{2m-1}\left(M+L(K-\x_\zs{*})\right)
+\frac{\beta}{2}x^{2m}\Chi_\zs{\{|x|\le K\}}\,-\,\frac{1}{2L}x^{2m}
\\[2mm]
&\le \A_\zs{m}-\frac{\beta}{2}x^{2m}\,,
\end{align*}
where 
$$
\A_\zs{m}=\left(2(\x_\zs{*}+ML)\right)^{2m-1}
\left(2M+\x_\zs{*}
\left(L+L^{-1}\right)+2L^{2}M
\right)
$$
From here it follows, that
$$
\dot{z}_\zs{t}(m)\le 2m\,\A_\zs{m}-L^{-1}\,
m\,z_\zs{t}(m)+m(2m-1)\sigma^{2}_\zs{\max}\,z_\zs{t}(m-1)\,.
$$
We can rewrite this inequality as follows
$$
\dot{z}_\zs{t}(m)=-L^{-1}mz_\zs{t}(m)\,+\,m(2m-1)\sigma_1^2z_\zs{t}(m-1)\,+\,\psi_t\,,
$$
where $\sup_\zs{t\ge 0} \psi_t\le 2m\,\A_\zs{m}$. This equality provides
\begin{align*}
z_\zs{t}(m)&=z_\zs{0}(m)e^{-mL^{-1}\,t}\,+\,
m(2m-1)\sigma^{2}_\zs{\max}\int^{t}_\zs{0}\,e^{-mL^{-1}(t-s)}z_\zs{s}(m-1)\d s\\[2mm]
&+\int^{t}_\zs{0}\,e^{-mL^{-1}(t-s)}\psi_\zs{s} \d s
\\[2mm]
&
\le m(2m-1)\sigma^{2}_\zs{\max}\int^{t}_\zs{0}\,e^{-mL^{-1}(t-s)}z_\zs{s}(m-1)\d s\,+
\B_\zs{m}\,,
\end{align*}
where 
$\B_\zs{m}=y^{2m}_\zs{0}+2\A_\zs{m}L$. Setting $\B_\zs{0}=1$ and
resolving this inequality by recurrence yields
$$
z_\zs{t}(m)\le 4\,(2m-1)!!\sum_\zs{j=0}^{m}\,
\left(\sigma^{2}_\zs{\max}L\right)^{m-j}\B_\zs{j}\,.
$$
It is easy to see, that
$$
\B_\zs{m}\le 4\left(\max\left(|y_0|^2, 4(\x_\zs{*}+ML)^{2}\right)\right)^{m}
\,.
$$
Therefore
$$
\sup_\zs{t\ge 0}z_\zs{t}(m)\le 4(m+1)(2m-1)!!\,\rho^{2m}
\le 4(2m)^{m}\,\rho^{2m}
\,,
$$
where $\rho$ is defined in \eqref{sec:Pa.5}.
Hence Proposition~\ref{Pr.sec:A.3}.
  \endproof

\medskip

\subsection{Properties of the function \eqref{sec:Pr.4}}

\begin{lemma}\label{Le.sec:A.1}
For any integrated function $\phi$ the solution \eqref{sec:Pr.4}
is uniform bounded, i.e.
$$
\sup_\zs{\vartheta\in\Theta}
\,\sup_\zs{y\in\bbr}|v_\zs{\vartheta}(y)|\le r\,,
$$
where the upper bound $r$ is introduced in \eqref{sec:Pa.4}.
\end{lemma}
\proof
Firstly we note, that for any  $\vartheta$ from $\Theta$
and any intergated $\bbr\to\bbr$ function $\phi$
$$
|\pi_\zs{\vartheta}(\phi)|\le q^{*}\,|\phi|_\zs{1}\,.
$$
Moreover, by the definition of the parameter $\beta_\zs{1}$ we get
$$
2\sup_\zs{|u|\le \x_\zs{*}}
|S_\zs{1}(u)|
\le \beta_\zs{1}\,.
$$
Therefore, for $0\le u\le \x_\zs{*}$ we can estimate the function 
$v_\zs{\vartheta}$ as
$$
|v_\zs{\vartheta}(u)|\le
\frac{2e^{\x_\zs{*}\beta_\zs{1}}}{\sigma^{2}_\zs{\min}}
\left(
\left(
1+q^{*}\x_\zs{*}
\right)
|\phi|_\zs{1}
+\I(\phi)
\right)
\,,
$$
where $\beta_\zs{1}$ is given in \eqref{sec:Pa.2} and
$$
\I(\phi)=\int_\zs{\x_\zs{*}}^{\infty}\,
\left(
|\phi(y)|
+
q^{*}\,|\phi|_\zs{1}
\right)\,
e^{2\int^{y}_\zs{\x_\zs{*}}\,S_\zs{1}(z)\d z}\,\d y\,.
$$
To estimate this term note that similarly to \eqref{sec:A.6} we can
obtain that for any  $y\ge a\ge \x_\zs{*}$
\begin{equation}\label{sec:A.7}
2\int^{y}_\zs{a}\,S_\zs{1}(z)\d z
\le 
\beta_\zs{1}(y-a)-\beta_\zs{2}(y-a)^{2}\,.
\end{equation}
Using this inequlity for $a=\x_\zs{*}$, we get
\begin{align*}
\I(\phi)
&\le
\int_\zs{\x_\zs{*}}^{\infty}\,
|\phi(y)|
e^{\beta_\zs{1}(y-\x_\zs{*})-\beta_\zs{2}(y-\x_\zs{*})^{2}}\,\d y
+
q^{*}\,|\phi|_\zs{1}
\int^{\infty}_\zs{0}\,
e^{\beta_\zs{1} z-\beta_\zs{2} z^{2}}\,\d z
\\[2mm]
&
\le 
|\phi|_\zs{1}\sup_\zs{z\ge 0}\,e^{\beta_\zs{1} z-\beta_\zs{2} z^{2}}
+
q^{*}\,|\phi|_\zs{1}
\int_\zs{0}^{\infty}\,
e^{\beta_\zs{1} z-\beta_\zs{2} z^{2}}\,\d z\\[2mm]
&\le 
|\phi|_\zs{1}
\left(
\upsilon_\zs{1}
+
q^{*}
\upsilon_\zs{2}
\right)\,,
\end{align*}
where the parameters $\upsilon_\zs{1}$ and $\upsilon_\zs{2}$
are introduced in \eqref{sec:Pa.2}. 
Therefore, taking into account the definition \eqref{sec:Pa.4}, the last inequality 
implies
\begin{equation}\label{sec:A.8}
\sup_\zs{\vartheta\in\Theta}\sup_\zs{0\le u\le \x_\zs{*}}|v_\zs{\vartheta}(u)|
\le
r\,.
\end{equation}
\noindent
If $u\ge \x_\zs{*}$, then 
through the inequality \eqref{sec:A.7}
 we estimate the function $v_\zs{\vartheta}(u)$ from above as
$$
\sup_\zs{\vartheta\in\Theta}\sup_\zs{u\ge \x_\zs{*}}|v_\zs{\vartheta}(u)|
\le
\frac{2|\phi|_\zs{1}}{\sigma^{2}_\zs{\min}}
\left(\upsilon_\zs{1}+q^*\upsilon_\zs{2}\right)
\le r\,.
$$
Let now $u\le 0$. Taking into account that
$$
\int_\zs{R}\,\frac{\wt{\phi}(y)}{\sigma^{2}(y)}\,
\exp\{2\int^{y}_\zs{0}\,S_\zs{1}(z)\,\d z\}\,\d y=0\,,
$$
we can represent the function $v_\zs{\vartheta}$ as
$$
v_\zs{\vartheta}(u)\,=\,2\int^{\infty}_\zs{|u|}\,\frac{\wt{\phi}(-y)}{\sigma^{2}(-y)}\,
e^{-2\int^{y}_\zs{|u|}\,S_\zs{1}(-z)\,\d z}\,\d y\,.
$$
Similarly to \eqref{sec:A.7}, one can check directly, that for any $y\ge a\ge \x_\zs{*}$
$$
-2\int^{y}_\zs{a}\,S_\zs{1}(-z)\,\d z
\le \beta_\zs{1}(y-a)-\beta_\zs{2}(y-a)^{2}\,.
$$
Therefore, by the same way as in the proof of
\eqref{sec:A.8}  we can estimate the function $v_\zs{\vartheta}(u)$
as
$$
\sup_\zs{\vartheta\in\Theta}\sup_\zs{u\le 0}|v_\zs{\vartheta}(u)|
\le\,r\,.
$$
Hence Lemma~\ref{Le.sec:A.1}.
\endproof

\end{document}